\documentclass[11pt,a4paper]{article}
\usepackage[T1]{fontenc}
\usepackage{lmodern,amsmath,amssymb,amsthm,microtype}
\usepackage[margin=24mm]{geometry}
\usepackage{tikz}
\usetikzlibrary{arrows.meta,positioning,calc}
\usepackage[hidelinks]{hyperref}
\usepackage{enumitem}
\newtheorem{theorem}{Theorem}

\newtheorem{corollary}[theorem]{Corollary}
\theoremstyle{definition}

\theoremstyle{remark}

\newcommand{\NP}{N^+}
\newcommand{\NM}{N^-}
\title{Highly connected uniquely $k$-linked digraphs}
\author{J\o{}rgen Bang-Jensen\thanks{
Department of Mathematics and Computer Science, University of Southern Denmark, Odense Denmark and School of Mathematics, Shandong University, Jinan China (email: jbj@imada.sdu.dk)}}
\date{}

\begin{document}
\maketitle

\begin{abstract}
Extending a result of Thomassen, we prove that for every choice of positive integers $K\geq k$ there exist a $K$-strong digraph which has two disjoint sets $S,T$ each of size  $k$ so that the only linkage from $S$ to $T$ in $D$ is the one whose existence follows from Menger's theorem.

\end{abstract}

\noindent{}{\bf Keywords:} disjoint paths; k-linkage; NP-complete; uniquely k-linked digraphs

\section{Introduction}
Notation and terminology is consistent with \cite{bang2009}. A digraph $D$ is {\bf $\mathbf{k}$-strong} if $D=(V,A)$ has at least $k+1$ vertices and $D-X$ is strongly connected for every $X\subset V$ with $|X|\leq k-1$. Every undirected graph
$G$ corresponds to the symmetric digraph $\stackrel{\leftrightarrow}{G}$ that we obtain by replacing each edge $u$ of $G$ by a pair of opposite arcs $uv,vu$. Clearly $\stackrel{\leftrightarrow}{G}$ is $k$-strong if and only if $G$ is $k$-connected. Hence our discussion below also covers undirected graphs.\\

Let $D$ be a digraph and let $S=\{s_1,s_2,\ldots{},s_k\}$ and 
$T=\{t_1,t_2,\ldots{},t_k\}$ be disjoint sets of $k$ vertices each. An {\bf $\mathbf{(S,T)}$-linkage} in $D$   is a collection of disjoint paths $P_1,P_2,\ldots{},P_k$ that start in $S$ and end in $T$. It realizes the permutation $\pi\in\mathfrak S_k$ if the  path starting in  $s_i$ ends
in $t_{\pi(i)}$. Different linkages may realize the same permutation. By Menger's theorem an $(S,T)$-linkage exists for every choice of disjoint $k$-sets $S,T$ when $D$ is $k$-strong. A digraph $D$ is {\bf $\mathbf{k}$-linked} if we can choose the paths $P_1,\ldots{},P_k$ above  such that $P_i$ starts in $s_i$ and ends in $t_{\tau(i)}$ for all $i\in [k]$ for each of the $k!$  permutations $\tau\in \mathfrak S_k$. At the other extreme, we call a $k$-strong digraph $D=(V,A)$ {\bf uniquely $\mathbf{k}$-linked} if there is some choice of disjoint $k$-sets $S,T\subset V$ such that  every $(S,T)$-linkage  in $D$ links $s_i$ to $t_{\pi(i)}$ for some fixed permutation $\pi\in\mathfrak S_k$. 

The {\bf $\mathbf{k}$-linkage problem} is as follows: decide for a given  (di)graph $G$ and distinct vertices $s_1,s_2,\ldots{},s_k$,$t_1,t_2,\ldots{},t_k$; whether $G$ has a collection of disjoint paths $P_1,P_2,\ldots{},P_k$ such that $P_i$ is an $(s_i,t_i)$-path.
Robertson and Seymour \cite{robertsonJCT63} proved that the $k$-linkage problem
is polynomial for all fixed $k$ for undirected graphs (symmetric digraphs).

Sufficient conditions for an undirected graph to be $k$-linked have been studied intensively. The best bound  is due to Thomas and Wollan \cite{thomasEJC26} who proved that  every $2k$-connected graph of minimum degree at least $10k$ is $k$-linked. 

For digraphs the situation is very different. Fortune, Hopcroft and Wyllie \cite{fortuneTCS10} proved that the 2-linkage problem is NP-complete for digraphs.
They also gave a polynomial algorithm for every fixed $k$ for the $k$-linkage problem for acyclic digraphs. An important result of Thomassen shows that for  digraphs there is no connectivity which is sufficient to guarantee that a digraph is  2-linked. In fact Thomassen proved that even more holds:

\begin{theorem}[Thomassen]\label{thm:CTnon2link}\cite{thomassenC11}
\begin{itemize}
   \item[(a)]  For every integer  $K\geq 1$ there exists a $K$-strong digraph which is not 2-linked.
   \item[(b)] For every integer $K\geq 1$ there exists a $K$-strong digraph $D$ with two vertices $u$ and $v$ such that $D$ has no directed cycle containing both $u$ and $v$.
   \item[(c)] For every integer $K\geq 1$ the 2-linkage problem is NP-complete for $K$-strong digraphs.
\end{itemize}
    
\end{theorem}

By (a) above, for every $K$ there exists a $K$-strong digraph $D$ so that, for some choice of terminal sets $S=\{s_1,s_2\}$ and $T=\{t_1,t_2\}$ the only $(S,T)$-linkage in $D$ links $s_i$ to $t_{3-i}$ for $i\in [2]$.  When $k=2$, there are only $2!=2$ different permutations of $T$ that an $(S,T)$-linkage may realize  so if $D$ is not 2-linked, then it is automatically uniquely 2-linked. As $k$ grows, the number of such possible permutations of a $k$-set grows very quickly, so it is tempting to think that there should exists an integer $k_0$ such that for all $k\geq k_0$, no  $k$-strong digraph is uniquely $k$-linked. However, as we show in this note, this is not the case. Building on Theorem \ref{thm:CTnon2link} (b) we 
construct arbitrarily highly connected uniquely-$k$-linked digraphs for all integers $k\geq 3$. Using Theorem \ref{thm:CTnon2link}(c) we also show that deciding whether an $r$-strong digraph not uniquely $k$-linked is NP-complete for all integers $r,k\geq 2$.

\section{Highly connected uniquely $\mathbf{k}$-linked digraphs }

\begin{theorem}\label{thm:main}
For every $r\ge1$ and $k\ge2$, there is a finite $r$-strong digraph $D$ with
disjoint vertex sets  $S=\{s_1,\ldots,s_k\}$, $T=\{t_1,\ldots,t_k\}$ such that an identity
$(S,T)$-linkage exists and every $(S,T)$-linkage realizes the identity
permutation. Consequently, $D$ is uniquely $k$-linked.
\end{theorem}

\begin{proof}

 We use the digraph guaranteed by Theorem \ref{thm:CTnon2link}(b) in our construction.

Fix $r,k$ and choose as the building block with $K=r+2$ a $K$-strong digraph $H$ with two vertices $u,v$ that do not lie on a common cycle. Clearly the vertices $u,v$ are nonadjacent since an arc between them, together with a
return path in $H$, would give such a cycle.\\

Take disjoint copies $H_1,\ldots,H_k$ of $H$ and let  $u_i,v_i$ denote the copies of $u,v$ in $H_i$ for $i\in [k]$.
Set
\[
R_i=H_i-\{u_i,v_i\},\qquad
A_i=\NP_{H_i}(v_i),\qquad B_i=\NM_{H_i}(v_i).
\]
As $u_i$ and $v_i$ are not adjacent, both $A_i$ and $B_i$ are subsets of $V(R_i)$; they need not be disjoint.
For each $i\in [k]$ we form a digraph  $D_i$ from $R_i$ by adding two nonadjacent vertices $s_i,t_i$.
For each $z\in\{s_i,t_i\}$, prescribe
\[
\NP_{D_i}(z)=\NP_{H_i}(u_i),\qquad
\NM_{D_i}(z)=\NM_{H_i}(u_i),
\]
and keep all the original arcs inside $R_i$.
Finally we obtain the digraph $D$ by taking  the disjoint union of  $D_1,\ldots{},D_k$ and adding  all arcs
\[
\boxed{\quad xy\quad\text{with }x\in B_i,\ y\in A_j,\ i\ne j.\quad}
\]
These are the only arcs between different $D_i$'s. Thus every arc leaving  $V(D_i)$ has its tail  in $B_i$
and every arc entering  $V(D_i)$ lie has its head in $A_i$.

\emph{Strong connectivity.}
Since $H_i$ is $(r+2)$-strong, each $R_i$ is $r$-strong and since $u_i,v_i$ are nonadjacent and are copies of $u,v$ in the $(r+2)$-strong digraph $H$, each of the new
terminals has at least $r+2$ in-neighbours and at least $r+2$ out-neighbours
in $R_i$. Therefore each $D_i$ is $r$-strong: after deleting fewer than $r$
vertices, its surviving core is strong and every surviving terminal still
has both an in-neighbour and an out-neighbour in that core.

Let $X\subseteq V(D)$ with $|X|<r$. Every $D_i-X$ is strong, and
\[
|A_i|\ge r+2,\qquad |B_i|\ge r+2.
\]
Thus $A_i-X$ and $B_i-X$ are nonempty. For every $i\ne j$, an arc remains
from $B_i-X$ to $A_j-X$. It follows that $D-X$ is strong, so $D$ is
$r$-strong.

\emph{Existence of an identity linkage.}
For each $i\in [k]$ choose an $(s_i,t_i)$-path $P_i$ entirely inside each strong piece $D_i$.
The pieces are vertex-disjoint, so the chosen paths $P_1,\ldots{},P_k$ form an $(S,T)$-linkage for the sets 
$S=\{s_1,s_2,\ldots{},s_k\}$ and $T=\{t_1,t_2,\ldots{},t_k\}$

\emph{Uniqueness of the terminal permutation.}
Suppose an $(S,T)$-linkage $Q_1,\ldots,Q_k$ realizes a permutation $\pi$,
where $Q_j$ starts at $s_j$ and ends at $t_{\pi(j)}$.
Assume that $\pi(i)\ne i$ for some $i$, and put $h=\pi^{-1}(i)$.
Then $h\ne i$, so $Q_i$ and $Q_h$ are distinct, vertex-disjoint paths.

The path $Q_i$ starts in $s_i\in V(D_i)$ and ends outside of $D_i$. Let $b\in B_i$ be the
vertex immediately before its \emph{first exit} from $D_i$. The prefix
\[
L=Q_i[s_i,b]
\]
is an $(s_i,b)$-path inside  $D_i$. Similarly, $Q_h$ starts outside $D_i$ and ends at $t_i$.
Let $a\in A_i$ be the vertex immediately after its \emph{last entry} into
$D_i$. The suffix
\[
M=Q_h[a,t_i]
\]
is an $(a,t_i)$-path which also lies in $D_i$.

Every terminal is an endpoint of one of the pairwise disjoint paths.
Hence none can occur internally on another path. In particular, all vertices
of $L-s_i$ and $M-t_i$ lie in $R_i$, and these two vertex sets are disjoint.

Write $L=s_i x_1\cdots x_p$, where $x_p=b$, and
$M=y_1\cdots y_qt_i$, where $y_1=a$. Here $p,q\ge1$.
In the original copy $H_i$, we obtain the directed paths
\[
u_i x_1\cdots x_p v_i
\qquad\text{and}\qquad
v_i y_1\cdots y_q u_i.
\]
The arcs incident with $u_i$ exist by the  definition of the in- and out-neighbourhoods of $s_i,t_i$;
the other two end arcs exist because $b\in\NM_{H_i}(v_i)$ and
$a\in\NP_{H_i}(v_i)$. Their interiors are disjoint, so their union is a
directed cycle through $u_i$ and $v_i$, a contradiction, since $H$ has no directed cycle through $u$ and $v$. Therefore we must have 
$\pi(i)=i$ for every $i$, proving that $D$ is uniquely $k$-linked.
\end{proof}

It is not difficult to see that Theorem \ref{thm:CTnon2link}(c) implies the following. We provide a proof for compleness.

\begin{theorem}\label{thm:NPCcycleKstrong}
For every fixed integer $r\geq 1$, it is NP-complete to decide whether
an $r$-strong digraph contains a directed cycle through two specified
distinct vertices.
\end{theorem}

\begin{proof}
Membership in NP follows by supplying the cycle.
For NP-hardness, apply Theorem \ref{thm:CTnon2link}(c) with $K=r+4$.
Let $(D;s_1,t_1,s_2,t_2)$ be an instance of the resulting restricted
2-linkage problem, and put
\[
 W=V(D)\setminus\{s_1,t_1,s_2,t_2\}.
\]
Construct a digraph $D'$ on $W\cup\{x,y\}$, where $x,y$ are new vertices.
Keep all arcs of $D[W]$, and add precisely the following arcs incident
with $x$ or $y$:
\[
\begin{array}{lll}
 xv &\text{if }s_1v\in A(D),&v\in W,\\
 vx &\text{if }vt_2\in A(D),&v\in W,\\
 yv &\text{if }s_2v\in A(D),&v\in W,\\
 vy &\text{if }vt_1\in A(D),&v\in W.
\end{array}
\]
In addition, add $xy$ if $s_1t_1\in A(D)$, and add $yx$ if
$s_2t_2\in A(D)$.
Thus $x$ combines the outgoing role of $s_1$ and the incoming role of
$t_2$, while $y$ combines the outgoing role of $s_2$ and the incoming
role of $t_1$. 

It is easy to check that $D'$ has a cycle through $x,y$ if and only if $D$ has disjoint paths $P_1,P_2$ such that $P_i$ is an $(s_i,t_i)$-path for $i\in [2]$.

Furthermore $D'[W]=D[W]$ is $r$-strong and since each of $x,y$ have at least $r+1$ in-and out-neighbours in $W$ it follows that $D'$ is $r$-strong (see also Proposition 5.8.5 in \cite{bang2009}).

The construction takes polynomial time. Together with the equivalence,
it gives the required polynomial reduction.
\end{proof}

Combining Theorem \ref{thm:NPCcycleKstrong} with the reduction in the proof of Theorem \ref{thm:main} we obtain the following.

\begin{corollary}
    For every fixed pair of positive integers $K\geq k$ is NP-complete to decide whether a $K$-strong digraph $D$ is uniquely $k$-linked.
\end{corollary}

\noindent{}{\bf Declaration about the use of AI:} AI was used in investigation of various properties of linkages in digraphs, in particular properties  of uniquely $k$-linked digraphs. This resulted in the reduction idea used in the proof of Theorem \ref{thm:main}. The author takes full responsibility for all parts of this paper.

\noindent{}{\bf Acknowledgement:} The author thanks Carsten Thomassen for stimulating discussions on linkages in digraphs which lead to the concept of uniquely $k$-linked digraphs.

\bibliographystyle{plain}

\begin{thebibliography}{1}

\bibitem{bang2009}
J.~Bang-Jensen and G.~Gutin.
\newblock {\em {Digraphs: Theory, Algorithms and Applications, 2nd Edition}}.
\newblock Springer-Verlag, London, 2009.

\bibitem{fortuneTCS10}
S.~Fortune, J.E. Hopcroft, and J.~Wyllie.
\newblock {The directed subgraph homeomorphism problem}.
\newblock {\em Theor. Comput. Sci.}, 10:111--121, 1980.

\bibitem{robertsonJCT63}
N.~Robertson and P.D. Seymour.
\newblock {Graph minors. XIII: The disjoint paths problem.}
\newblock {\em J. Combin. Theory Ser. B}, 63:65--110, 1995.

\bibitem{thomasEJC26}
R.~Thomas and P.~Wollan.
\newblock An improved linear edge bound for graph linkages.
\newblock {\em Eur. J. Comb.}, 26(3-4):309--324, 2005.

\bibitem{thomassenC11}
C.~Thomassen.
\newblock {Highly connected non-{$2$}-linked digraphs}.
\newblock {\em Combinatorica}, 11(4):393--395, 1991.

\end{thebibliography}

\end{document}